\documentclass[12pt]{amsart}

\usepackage{amsmath,amssymb,amsthm}

\usepackage{amsmath, amssymb}
\usepackage{amsthm, amsfonts, mathrsfs}
\usepackage{amsfonts,graphicx}
\theoremstyle{plain}
\newtheorem{Theo}{Theorem}[section]
\newtheorem{lem}[Theo]{Lemma}
\newtheorem{cor}[Theo]{Corollary}
\newtheorem{prop}[Theo]{Proposition}

\theoremstyle{plain}
\theoremstyle{definition}

\theoremstyle{remark}
\newtheorem{Rema}[Theo]{Remark}
\newtheorem*{rema*}{Remark}

\newcommand{\NN}{\mathbb{N}}

\newcommand{\RR}{\mathbb{R}}

\author[H. Abidi]{Hammadi Abidi}
\address{Facult\'e des Sciences  de Tunis}
\email{Hamadi.Abidi@fst.rnu.tn}
\author[T. Hmidi]{Taoufik Hmidi}
\address{IRMAR, Universit\'e de Rennes 1\\ Campus de
Beaulieu\\ 35~042 Rennes cedex\\ France}
\email{thmidi@univ-rennes1.fr}
\author[S. Keraani]{Sahbi Keraani}
\address{IRMAR, Universit\'e de Rennes 1\\ Campus de
Beaulieu\\ 35~042 Rennes cedex\\ France}
\email{sahbi.keraani@univ-rennes1.fr}

\date{}

\begin{document}

\title[On the axisymmetric Boussinesq system]
{On the global regularity of axisymmetric Navier-Stokes-Boussinesq system}
\maketitle
\begin{abstract}
In this paper we prove a global well-posedness result  for   tridimensional  Navier-Stokes-Boussinesq system  with axisymmetric initial data. This system couples Navier-Stokes equations with a transport equation governing the density. 
\end{abstract}
\tableofcontents
%%%%%%%%%%%%%%
\section{Introduction}

The purpose of this paper is to study the global well-posedness  for three-dimensional Boussinesq system
in the whole space  with axisymmetric initial data. This system is described as follows,
\begin{equation}
      \label{bs}
\left\{ \begin{array}{ll}
\partial_{t}v+v\cdot\nabla v-\Delta v+\nabla p=\rho e_z,\quad (t,x)\in \RR_+\times\RR^3,\\
\partial_{t}\rho+v\cdot\nabla\rho=0,\\
\textnormal{div}\,  v=0,\\
v_{| t=0}=v_{0},   \quad \rho_{| t=0}=\rho_{0}.
\end{array} \right.
\end{equation}
Here, the velocity  $v=(v^1,v^2,v^3)$ is a three-component vector field with zero divergence. The scalar
function $\rho$ denotes the density which is transported by the flow and acting for  the first equation of \eqref{bs}  only in the vertical direction given by $e_z$. The pressure $p$ is a scalar function related to the
unknowns $v$ and $\rho$ through an elliptic equation.
Remark that when the initial density $\rho_0$ is identically zero  then the above system is reduced to the
classical Navier-Stokes system  which is widely studied:
\begin{equation}
      \label{ns}
\left\{ \begin{array}{ll}
\partial_{t}v+v\cdot\nabla v-\Delta v+\nabla p=0\\
\textnormal{div}\,  v=0\\
v_{| t=0}=v_{0}.
\end{array} \right.
\end{equation}
 Recall that the existence of global weak solutions in the energy space for \eqref{ns}
 goes back to J. Leray \cite{Leray} in the last century. However the uniqueness of these solutions
 is only known in space dimension two. It is also well-known that  smooth solutions are global
 in dimension two and  for higher dimensions when the data are small in some critical spaces; see for instance \cite{Lemar} for more detailed  discussions. Although
 the breakdown of smooth solutions with large initial data is still now an open problem
 some partial results are known outside the context of small data. We refer for instance to
  recent papers of J.-Y. Chemin and I. Gallagher \cite{CG,CG1} where  global existence in dimension
 three is established for special structure of initial data which are not small in any critical space.

\
There is  an interesting   case  of    global existence for \eqref{ns}  corresponding  to large initial data but with special
geometry, called axisymmetric without swirl.
Before going further in the details  let us give some general statements  about Navier-Stokes system
in space dimension three. First, we start with introducing the vorticity which is a physical
quantity that plays a significant role in the theory of global existence;
 for a given vector field $v$ the vorticity  $\omega$ is  the vector defined by
 $\omega={\rm curl}\,  v.$ Thus we get from \eqref{ns}  the vorticity  equation
$$
\partial_t \omega+v\cdot\nabla\omega-\Delta\omega=\omega \cdot\nabla v.
$$
According to Beale-Kato-Majda criterion the formation of
singularities in finite time is due to the accumulation of the
vorticity. In other words, to have global existence it suffices to
bound for every time the quantity $\|\omega(t)\|_{L^\infty}$.
However the main difficulty arising in dimension  three is the lack
of information about the manner that the vortex-stretching term
$\omega \cdot\nabla\,v$ affects the dynamic of the fluid.

\ For the geometry of axisymmetric flows without swirl  we have a 
cancellation in the stretching term giving rise to  new conservation
laws. We say that a  vector field $v$ is  axisymmetric  if it has
the form:
$$
v(t,x) = v^r(t,r, z)e_r + v^z (t,r, z)e_z,\quad x=(x_{1},x_{2},z),\quad r=({x_{1}^2+x_{2}^2})^{\frac12},
$$
where  $\big(e_r, e_{\theta} , e_z\big)$ is the cylindrical basis of
$\mathbb R^3$ and the components $v^r$ and $v^z$ do not depend on
the angular variable.   The main feature  of axisymmetric  flows
arises in the vorticity  which takes the form,
$$
\omega=(\partial_zv^r-\partial_rv^z) e_{\theta}:=\omega_\theta e_\theta
$$
and satisfies
  \begin{equation*}
\partial_t \omega +v\cdot\nabla\omega-\Delta\omega =\frac{v^r}{r}\omega.
\end{equation*}
Since the Laplacian operator has the form $\Delta=\partial_{rr}+\frac1r\partial_r+\partial_{zz}$
in the cylindrical coordinates then the  component $\omega_\theta$ of the vorticity will satisfy
\begin{equation}
 \label{tourbillon}
\partial_t \omega_\theta +v\cdot\nabla\omega_\theta-\Delta\omega_\theta
+\frac{\omega_\theta}{r^2} =\frac{v^r}{r}\omega_\theta.
\end{equation}
Consequently, the quantity $\Gamma:=\frac{\omega_\theta}{r}$  obeys to the equation
\begin{equation*}
%\label{equation_importante}
\partial_t\Gamma+v\cdot\nabla\Gamma-\Delta\Gamma-\frac2r\partial_r\Gamma=0.
\end{equation*}
Obviously, we have for $p\geq1,$  ${\int\frac2r\partial_r\Gamma |\Gamma|^{p-1}(\textnormal{sign } \Gamma) dx\leq0}$
and  then we deduce that for all $p\in[1,\infty]$
$$
\|\Gamma(t)\|_{L^p}\le\|\Gamma_0\|_{L^p}.
$$
It was shown by M. Ukhoviskii and V. Yudovich \cite{Ukhovskii}  and independently by O. A. Ladyzhenskaya \cite{LA}  that these  conservation laws are  strong enough to prevent the formation of singularities in
finite time for axisymmetric flows. More precisely, the system \eqref{ns} has a unique global solution for  $v_0\in H^1$ such that
$\omega_0,\frac{\omega_0}{r}\in L^2\cap L^\infty.$  Remark that in term of Sobolev
regularities these assumptions are satisfied when $v_0\in H^s$ with $s>\frac72.$
Few decades afterwards,   S. Leonardi, J. M\`alek, J. Nec\u as and M. Pokorn\'y  \cite{lmnp}
weakened the initial regularity for   $v_0\in H^2$. This result was recently improved
by \mbox{H. Abidi \cite{a}}  for  $v_0\in H^{\frac{1}{2}}$.

Let us now come back to our first problematic which is the study of
global well-posedness for the Boussinesq system \eqref{bs}. In space
dimension two many papers are recently devoted to this problem and
the study seems to be in a satisfactory state. More precisely we
have global existence  in different function spaces  and for
different viscosities, we refer  for example to \cite{ah,Brenier,
ha,dp1,dp,dp2,hk1,hk,HKR1,HKR2}.

In the case of space dimension three few results are known about global existence. We recall
the result of R. Danchin and M. Paicu \cite{dp}. They proved a global well-posedness result for
small initial data belonging  to some critical Lorentz spaces.
% \mbox{namely that $\|v_0\|_{L^{3,\infty}}+\|\rho_0\|_{L^1}<<1$.}

  Our goal here is to study the global existence for the system \eqref{bs} with
  axisymmetric  initial data, which means that the velocity $v_0$ is assumed to be an axisymmetric vector
  field without swirl  and the density $\rho_0$ depends only on $(r,z)$. It should be mentioned that
  this structure is preserved for strong solutions in their  lifespan. Before stating our
  main result we denote by $\Pi_z$ the orthogonal projector over the axis $(Oz).$ Our
  result reads as follows.
 \begin{Theo}\label{thm1}
 Let $v_0\in H^1$ be an axisymmetric vector field with zero divergence and such that $\frac{{\omega_0}}{r}\in L^2.$
 Let $\rho_0\in L^2\cap L^\infty$ depending only on $(r,z)$ and such that $\hbox{supp }\rho_0$
 does not intersect the axis $(Oz)$ and $\Pi_z(\hbox{supp }\rho_0)$ is a compact set.
 Then the system \eqref{bs} has a unique global solution $(v,\rho)$ such that
 $$
 v\in{C}(\RR_+;H^1)\cap L^1_{\textnormal{loc }}(\RR_+;W^{1,\infty}),$$
 $$ \frac{\omega}{r}\in L^\infty_{\textnormal{loc}}(\RR_+;L^2),\quad
 \rho\in L^\infty_{\textnormal{loc}}(\RR_+;L^2\cap L^\infty).
 $$
  \end{Theo}
   The assumption that the density $\rho$ is zero in some region of the space $\RR^3$
   is not very meaningful from a physical  point  of view. However we can relax this
   hypothesis and extend our result to  more general case:  the initial density
   is assumed to be constant near the axis $(Oz)$ and for large value \mbox{of $z$.} More precisely, we have the following result.
  \begin{cor}
   Let $v_0\in H^1$ be an axisymmetric vector field with zero divergence and such that $\frac{{\omega_0}}{r}\in L^2.$
   Let $\rho_0\in L^\infty$   depending only on $(r,z)$ such that $\rho_0\equiv c_0$, for some
   constant $c_0$  in a region of type $\big\{x; r\le r_0, |z|\geq |z_0|  \big\},$ with $r_0>0$ and
   $\rho_0-c_0\in L^2$. Then the system \eqref{bs} has a unique global solution $(v,\rho)$ such that
 $$
 v\in{C}(\RR_+;H^1)\cap L^1_{\textnormal{loc }}(\RR_+;W^{1,\infty}),$$
 $$ \frac{\omega}{r}\in L^\infty_{\textnormal{loc}}(\RR_+;L^2),\quad
 \rho-c_0\in L^\infty_{\textnormal{loc}}(\RR_+;L^2\cap L^\infty).
$$
  \end{cor}
  The proof of this corollary is an immediate  consequence of Theorem \ref{thm1}. Indeed, we set
  $\bar\rho(t,x)=\rho(t,x)-c_0,$ then the system \eqref{bs} is reduced to
  \begin{equation*}
\left\{ \begin{array}{ll}
\partial_{t}v+v\cdot\nabla v-\Delta v+\nabla p-c_0e_z=\bar\rho e_z\\
\partial_{t}\bar\rho+v\cdot\nabla\bar\rho=0\\
\textnormal{div}\,  v=0\\
v_{| t=0}=v_{0},   \quad \bar\rho_{| t=0}=\rho_{0}-c_0.
\end{array} \right.
\end{equation*}
Now by  changing  the pressure $p$ to $\bar p=p-c_0 z$  we get the same system \eqref{bs}
and therefore we can apply the results of Theorem \ref{thm1}.

  Next we shall briefly discuss the new  difficulties  that one should  deal with compared to
  the system \eqref{ns}.
  First we start with writing  the analogous equation to  \eqref{tourbillon} for the vorticity.
  An easy computation gives
  $$
  \textnormal{curl}(\rho e_z)=\begin{pmatrix}
      \partial_2\rho   \\
     \partial_1\rho\\
     0
\end{pmatrix}=-(\partial_r\rho) e_\theta.
 $$
This yields to
\begin{equation}
 \label{tourbillon1}
\partial_t \omega_\theta +v\cdot\nabla\omega_\theta-\Delta\omega_\theta
+\frac{\omega_\theta}{r^2} =\frac{v^r}{r}\omega_\theta-\partial_r\rho.
\end{equation}
It follows that  the evolution of the quantity $\Gamma:=\frac{\omega_\theta}{r}$
is governed by  the equation
\begin{equation}
\label{equation_i}
\partial_t\Gamma+v\cdot\nabla\Gamma-\Delta\Gamma-\frac2r\partial_r\Gamma=-\frac{\partial_r\rho}{r}\cdot
\end{equation}
 Since the density $\rho$ satisfies a  transport equation then  the only conserved quantities
 that one should  use  are $\|\rho(t)\|_{L^p}$ for every $p\in[1,\infty].$  Loosely speaking,
 the source term  $\frac{\partial_r\rho}{r}$  in  the equation \eqref{equation_i} that one need to estimate has the scale of  $\Delta\rho$ and
 then to estimate $\Gamma$ in some Lebesgue spaces one can try for example  the maximal regularity of the
 heat semigroup. However it is not at all clear whether we can   prove a suitable maximal regularity because the involved
 elliptic operator $-\Delta-\frac{2}{r}\partial_r$ has singular \mbox{coefficients . }
Now, by
 taking the $L^2$-inner product of \eqref{equation_i} with $\Gamma$ we are led to estimate the quantity
 $\|\rho/r\|_{L^2}$ which has the scaling of $\|\nabla\rho\|_{L^2}$. In other words we have
 $$
 \|\Gamma(t)\|_{L^2}^2+\int_0^t\|\Gamma(\tau)\|_{\dot H^1}^2d\tau\le \|\Gamma_0\|_{L^2}^2
 +\int_0^t\|(\rho/r)(\tau)\|_{L^2}^2d\tau.
 $$
 Remark that the quantity $\|\rho/r\|_{L^2}$ is not well-defined if we don't assume at least
 that $\rho(t,0,z)=0.$ As we have seen previously   in the statement of Theorem \ref{thm1} we need more
 than this latter condition: the initial density   $\rho_0$ has to be supported far from the
 axis $(Oz).$  To follow the evolution of the quantity $\|\rho(t)/r\|_{L^2}$ the approach that consists 
 in writing the equation of $\frac{\rho(t)}{r}$ fails because it gives an exponential growth,
 $$
 \|\rho(t)/r\|_{L^2}\le  \|\rho_0/r\|_{L^2} e^{\|v^r/r\|_{L^1_tL^\infty}}.
 $$
 \
 Our idea to perform this growth relies in studying  some  dynamic aspects   of the support of $\rho(t)$ which  is nothing but the transported of the initial support by the flow.
   We need particularly to  get a proper low bound about   the distance from
   the axis $(Oz)$ of  the support of $\rho(t),$ see Proposition \ref{prop2}. To reach this target  and once again the axisymmetry property of the flow will   play a
   crucial role  since the trajectories of the particles are contained in the meridional plane.
   This approach allows us to improve the exponential growth for $ \|\rho(t)/r\|_{L^2}$ to a
   quadratic estimate. Roughly speaking we obtain the estimate
   $$
    \|\rho(t)/r\|_{L^2}\le C_0\|v^r/r\|_{L^1_tL^\infty}\big(1+\|v\|_{L^1_tL^\infty}\big).
   $$
   \
  The rest of the paper is organized as follows. In section $2$ we  recall some basic ingredients
  of Littlwood-Paley theory. In section $3$ we study some qualitative and analytic   properties
  of the flow associated to an axisymmetric vector field. In \mbox{section $4$} we give some global
  a priori estimates. The proof of Theorem \ref{thm1} is done in Section $5.$
  \section{Preliminaries}
  Throughout this paper, $C$ stands for some real positive constant which may be different
  in each occurrence and $C_0$ for a positive constant depending on the initial data.   We shall
  sometimes alternatively use the notation $X\lesssim Y$ for an inequality of type $X\leq CY$.

 Let us start with  a classical dyadic decomposition of the whole space (see  \cite{che1}):
there exist two radial  functions  $\chi\in \mathcal{D}(\mathbb R^3)$ and
$\varphi\in\mathcal{D}(\mathbb R^3\backslash{\{0\}})$ such that
\begin{enumerate}
\item
$\displaystyle{\chi(\xi)
+\sum_{q\geq0}\varphi(2^{-q}\xi)=1}$\;
$\forall\xi\in\mathbb R^3,$
\item
$ \textnormal{supp }\varphi(2^{-p}\cdot)\cap
\textnormal{supp }\varphi(2^{-q}\cdot)=\varnothing,$ if  $|p-q|\geq 2$,\\

\item
$\displaystyle{q\geq1\Rightarrow \textnormal{supp}\chi\cap
\textnormal{supp }\varphi(2^{-q})=\varnothing}$.
\end{enumerate}
For every $u\in{\mathcal S}'(\mathbb R^3)$ we define the nonhomogeneous Littlewood-Paley operators by,
$$
\Delta_{-1}u=\chi(\hbox{D})u;\, \forall
q\in\mathbb N,\;\Delta_qu=\varphi(2^{-q}\hbox{D})u\; \quad\hbox{and}\quad
S_qu=\sum_{-1\leq j\leq q-1}\Delta_{j}u.
$$
One can easily prove that for every tempered distribution $u,$
\begin{equation}\label{dr2}
u=\sum_{q\geq -1}\Delta_q\,u.
\end{equation}

In the sequel we will  make an extensive use of Bernstein inequalities (see for \mbox{example \cite{che1}}).
\begin{lem}\label{lb}\;
 There exists a constant $C$ such that for $k\in\NN$, \mbox{$1\leq a\leq b$}   and $u\in L^a$, we have
\begin{eqnarray*}
\sup_{|\alpha|=k}\|\partial ^{\alpha}S_{q}u\|_{L^b}
&\leq&
C^k\,2^{q(k+3(\frac{1}{a}-\frac{1}{b}))}\|S_{q}u\|_{L^a},
\end{eqnarray*}
and for $q\in\NN$
\begin{eqnarray*}
\ C^{-k}2^
{qk}\|{\Delta}_{q}u\|_{L^a}
&\leq&
\sup_{|\alpha|=k}\|\partial ^{\alpha}{\Delta}_{q}u\|_{L^a}
\leq
C^k2^{qk}\|{\Delta}_{q} u\|_{L^a}.
\end{eqnarray*}
\end{lem}

Let us now introduce the basic tool of  the paradifferential calculus which is  Bony's
decomposition  \cite{b}.   It distinguishes in   a product
$uv$  three parts as follows:
$$
uv=T_u v+T_v u+\mathcal{R}(u,v),
$$
where
\begin{eqnarray*}
T_u v=\sum_{q}S_{q-1}u\Delta_q v, \quad\hbox{and}\quad \mathcal{R}(u,v)=
\sum_{q}\Delta_qu \widetilde \Delta_{q}v,
\end{eqnarray*}
$$
\textnormal{with}\quad {\widetilde \Delta}_{q}=\sum_{i=-1}^{1}\Delta_{q+i}.
$$
$T_{u}v$ is called paraproduct of $v$ by $u$ and  $\mathcal{R}(u,v)$ the remainder term.

 Let $(p,r)\in[1,+\infty]^2$ and $s\in\mathbb R,$ then the nonhomogeneous  Besov
\mbox{space $B_{p,r}^s$} is
the set of tempered distributions $u$ such that
$$
\|u\|_{B_{p,r}^s}:=\Big( 2^{qs}
\|\Delta_q u\|_{L^{p}}\Big)_{\ell^{r}}<+\infty.
$$
We remark that the usual Sobolev space $H^s$ agrees with  Besov space  $B_{2,2}^s$.
Also, by using the Bernstein inequalities  we get easily
$$
B^s_{p_1,r_1}\hookrightarrow
B^{s+3({1\over p_2}-{1\over p_1})}_{p_2,r_2}, \qquad p_1\leq p_2\quad and \quad  r_1\leq r_2.
$$

\

%Let us now define Lorentz spaces $L^{p,q}$. For a measurable function $f$ we define its nonincreasing rearrangement  by
%$$
%f^{\ast}(t):= \inf\Big\{s,\;\mu\big(\{ x,\; |f (x)| > s\}\big)\leq t\Big\},
%$$
%where $\mu$ denotes the usual Lebesgue measure.
%For $(p,q)\in [1,+\infty]^2,$ the Lorentz space $L^{p,q}$ is the set of  functions $f$
%such that $\|f\|_{L^{p,q}}<\infty,$ with
%$$
%\|f\|_{L^{p,q}}:=\left\lbrace
%\begin{array}{l}\displaystyle
%\Big(\int_0^\infty[t^{1\over p}f^{\ast}(t)]^q{dt\over t}\Big)^{1\over q}, \quad \hbox{for}\;1\leq q<\infty\\
%\displaystyle\sup_{t>0}t^{1\over p}f^{\ast}(t),\quad \hbox{for}\;q=\infty.
%\end{array}
%\right.
%$$
%Notice that we can also define Lorentz spaces  by real interpolation from Lebesgue spaces:
%$$
%(L^{p_0},L^{p_1})_{(\theta,q)}=L^{p,q},
%$$
%where
%$
%1\leq p_0<p<p_1\leq\infty,
%$
%$\theta$ satisfies ${1\over p}={1-\theta\over p_0}
%+{\theta\over p_1}$ and $1\leq q\leq\infty$.
%We have the classical properties:
%\begin{equation}\label{imbed0}
%\|uv\|_{L^{p,q}}
%\le
%\|u\|_{L^\infty}\|v\|_{L^{p,q}}.
%\end{equation}
%\begin{equation}
%\label{imbed23}L^{p,q}\hookrightarrow L^{p,q'},\forall\, 1\leq p\leq\infty; 1\leq q\leq q'\leq \infty\quad \hbox{and}\quad L^{p,p}=L^p.
%\end{equation}

%

%

\section{Study of the flow map}
The main goal of this section is to study
some geometric and analytic properties of the  generalized flow map associated to  an axisymmetric vector field.
$$
\psi(t,s,x)=x+\int_s^tv(\tau,\psi(\tau,s,x))d\tau.
$$
This part will be the cornerstone of the proof of Theorem \ref{thm1}.
We need first to recall some basic results about the generalized flow. If the vector field
$v$  belongs to $L^1_{\textnormal{loc}}(\RR, C^{1}_b)$ then the generalized flow is
uniquely determined and exists globally in time. Here $C^1_b$ denotes the space of functions
with continuous bounded  gradient. In addition, for every $t,s\in\RR$, $\psi(t,s)$ is a
diffeomorphism that preserves Lebesgue measure when   $\textnormal{div } v=0$  and $$
\psi^{-1}(t,s,x)=\psi(s,t,x).
$$
Now we define the distance from a given point $x$ to a subset $A\subset\RR^3$ by
$$
d(x,A):=\inf_{y\in A}\|x-y\|,
$$
where $\|\cdot\|$ is the usual Euclidian norm.  The distance between two subsets $A$ and
$B$ of $\RR^3$ is defined by
$$
d(A,B):=\inf_{x\in A, y\in B}\|x-y\|.
$$
The diameter of a bounded subset $A\subset \RR^3$ is defined by
$$
\textnormal{diam }A=\sup_{x,y\in A}\|x-y\|.
$$
Our  first result is the following.
\begin{prop}\label{prop1}
Let $v$ be a smooth axisymmetric vector field and $\psi(t,s)$ its flow.
Let $x\notin (Oz)$ and $r(x):=d(x,(Oz))$, then
\begin{enumerate}
\item For every $s\in\RR$, the trajectory $\Gamma_{x,s}:=\big\{\psi(t,s,x),t\in\RR\big\}$ is a smooth
curve contained in the meridional plan.

\item For every $s\in\RR$, the trajectory $\Gamma_{x,s}$ does not intersect the axis $(Oz)$, that is
$\Gamma_{x,s}\cap (Oz)=\varnothing.$ More precisely,
$$
r(x)\displaystyle{e^{-|\int_s^t\|\frac{v^r}{r}(\tau)\|_{L^\infty}d\tau|}}
\le
d(\psi(t,s,x),(Oz))\le r(x)e^{|\int_s^t\|\frac{v^r}{r}(\tau)\|_{L^\infty}d\tau|}.
$$
\end{enumerate}
\end{prop}
\begin{proof}
{\bf (1)} We start with the decomposition of  the vector $x$ in the cylindrical basis:
$$
x=r(x)\begin{pmatrix}
      \cos\theta_x   \\
     \sin\theta_x\\
     0
\end{pmatrix}+z_xe_z,\quad\hbox{with}\quad r_x>0.
$$
 We decompose also the generalized flow  $\psi(t,s,x)$ in the cylindrical coordinates
$$
\psi(t,s,x)=r(t,s,x)\begin{pmatrix}
      \cos(\theta(t,s,x))   \\
     \sin(\theta(t,s,x))\\
     0
\end{pmatrix}+z(t,s,x)e_z.
$$
Since $v$ is smooth in space-time variables then  the flow map $\psi(t,s,x)$  is also smooth. Now we intend to prove that $r(t,s,x)$ remains strictly positive for all $t,s\in\RR$. Assume that there exist $t_1,s_1\in\RR$ such that $\psi(t_1,s_1,x)$ belongs to the axis $(Oz).$ As the restriction of the vector field $v$  on the axis $(Oz)$ satisfies $v(0,0,z)=v^z(t,0,z) e_z$ then trajectory for every $x_0\in (Oz)$ lies in this same   axis
$$
\psi(t,s,x_0)=x_0+\big(0,0,\int_s^t v^z(\tau,s,x_0)d\tau\big).
$$
Clearly one can choose  $x_0$  such that $\psi(t_1,s_1,x_0)=\psi(t_1,s_1,x)$ and this contradicts the fact that the flow map is an homeomorphism. Consequently and from the smoothness  of the generalized flow map
one can prove easily that   the functions $r(t,s,x), \theta(t,s,x)$  and $z(t,s,x)$  are also smooth in each  variable. This allows us to justify the following computation
\begin{eqnarray*}
\partial_t\psi(t,s,x)&=&\partial_t r(t,s,x)\begin{pmatrix}
      \cos(\theta(t,s,x))   \\
     \sin(\theta(t,s,x))\\
     0
\end{pmatrix}\\
&+&r(t,s,x)\partial_t\theta(t,s,x)\begin{pmatrix}
     - \sin(\theta(t,s,x))   \\
     \cos(\theta(t,s,x))\\
     0
\end{pmatrix}
+\partial_tz(t,s,x)e_z.
\end{eqnarray*}
Since the vector field $v$ is axisymmetric without swirl then
$$
v(t,\psi(t,s,x))=v^r(t,r(t,s,x),z(t,s,x))\begin{pmatrix}
 \cos(\theta(t,s))   \\
  \sin(\theta(t,s))\\
     0
\end{pmatrix}+v^z(t,r(t,s,x),z(t,s,x))e_z.
$$

Thus we get by identification
\begin{eqnarray}\label{eq1}
\nonumber&&\partial_t r(t,s,x)=v^r(t,r(t,s,x),z(t,s,x)),\\
 && r(t,s,x)\partial_t\theta(t,s)=0,\\
\nonumber&&\partial_t z(t,s,x)=v^z(t,r(t,s,x),z(t,s,x)).
\end{eqnarray}
From the above discussion we have $r(t,s,x)>0$ for all $t,s\in\RR$ and therefore we get
$$
\theta(t,s,x)=\theta(s,s,x)=\theta_x,\,\forall t,s\in\RR.
$$
It follows that for every $s\in\RR$  the trajectory $\Gamma_{x,s}:=\big\{\psi(t,s,x),t\in\RR\big\}$ lies in the meridional plan.

{\bf{(2)}}
From the first equation of \eqref{eq1} we get
$$
r(t,s,x)=r(x)+\int_{s}^tr(\tau,s,x)^{-1}v^r\big(\tau,r(\tau,s),z(\tau,s)\big)r(\tau,s,x)d\tau.
$$
Using Gronwall lemma we get
$$
r(t,s,x)\le r(x)e^{|\int_s^t\|\frac{v^r}{r}(\tau)\|_{L^\infty}d\tau|}.
$$
This gives the second inequality of the r.h.s since $d(\psi(t,s,x),(Oz))=r(t,s,x).$
Now we apply the above inequality by taking $ \psi(s,t,x)$ instead of  $x$. Since
$\psi(t,s,\psi(s,t,x))=x$ then
$$
r(x)\le r(s,t,x)e^{|\int_s^t\|\frac{v^r}{r}(\tau)\|_{L^\infty}d\tau|}.
$$
Interchanging  $(s,t)$ and $(t,s)$ gives
$$
r(x)e^{-|\int_s^t\|\frac{v^r}{r}(\tau)\|_{L^\infty}d\tau|}\le r(t,s,x).
$$
This achieves the proof of Proposition \ref{prop1}.
\end{proof}
 Proposition \ref{prop1} will be of much use in deriving some a priori estimates about solutions of transport equations. The first application is given below,
\begin{prop}\label{prop2}
Let $v$ be a smooth axisymmetric vector field and $\rho$ a solution of the transport equation
\begin{displaymath}
\left\{ \begin{array}{ll}
\partial_{t}\rho+v\cdot\nabla \rho=0&\\
%\textnormal{div }v_\nu=0 &\\
{\rho}_{| t=0}=\rho_{0}.
\end{array} \right.
\end{displaymath}
\begin{enumerate}
\item
Assume that $d(\textnormal{supp }\rho_0,(Oz))=r_0>0.$ Then we have for every $t\geq0$
$$
 d(\textnormal{supp }\rho(t),(Oz))\geq\, r_0\ e^{-\int_0^t\|\frac{v^r}{r}(\tau)\|_{L^\infty}d\tau}.
$$
\item Denote by $\Pi_z$ the orthogonal projector over the axis $(Oz)$. We assume that
$\Pi_z(\textnormal{supp }\rho_0)$ is compact set with diameter $d_0.$ Then for every
$t\geq0$, $\Pi_z(\textnormal{supp }\rho(t))$ is compact set with diameter $d(t)$ such that
$$
d(t)\le d_0+2\int_0^t\|v(\tau\|_{L^\infty}d\tau.
$$

\end{enumerate}

\end{prop}
\begin{proof}
{\bf (1)} The solution $\rho$ of the transport equation is completely described
through the usual  flow $\psi$, that is, $\rho(t,x)=\rho_0(\psi^{-1}(t,x)).$  Here
$\psi(t,x):=\psi(t,0,x)$ where $\psi(t,s,x)$ is the generalized flow introduced in the
beginning  of this section. Thus it follows that
$\textnormal{supp }\rho(t)=\psi(t,\textnormal{supp }\rho_0).$
Let now $y\in \textnormal{supp }\rho(t)$ then by definition
\mbox{ $y=\psi(t,x)$} with $x\in \textnormal{supp }\rho_0$.
Therefore  $d(y,(Oz))=r(t,x)$  with $r(t,x):=r(t,0,x)$ and so from \mbox{Proposition \ref{prop1}} we get
\begin{eqnarray*}
d(y,(Oz))&\geq& d(x,(Oz))\,e^{-\int_0^t\|\frac{v^r}{r}(\tau)\|_{L^\infty}d\tau}\\
&\geq& d(\textnormal{supp }\rho_0,(Oz))\,e^{-\int_0^t\|\frac{v^r}{r}(\tau)\|_{L^\infty}d\tau}\\
&\geq&r_0 \,e^{-\int_0^t\|\frac{v^r}{r}(\tau)\|_{L^\infty}d\tau}.
\end{eqnarray*}
This concludes the first part.

{\bf(2)} Let $x,\tilde x\in \textnormal{supp }\rho_0$ and  denote by
$y(t)=\psi(t,x)$ and $\tilde y(t)=\psi(t,\tilde x)$ . We set successively
$z(t)$ and $\tilde z(t)$ the last component of $y(t)$ and $\tilde y(t).$
Using the equation \eqref{eq1} with $s=0$ we get
$$
 \dot z(t)=v^z(t,r(t,x),z(t)).
$$
Integrating this differential equation we get
$$
z(t)=z(0)+\int_0^t v^z(\tau,r(\tau,x),z(\tau))d\tau.
$$
It follows that
$$
|z(t)-\tilde z(t)|\le|z(0)-\tilde z(0)|+2\int_0^t\|v(\tau)\|_{L^\infty}d\tau.
$$
This yields
$$
\textnormal{diam}(\Pi_z(\textnormal{supp }\rho(t)))\le\textnormal{diam}(\Pi_z(\textnormal{supp }\rho_0))+2\int_0^t\|v(\tau)\|_{L^\infty}d\tau.
$$
The proof is now completed.
\end{proof}
Now we are in a position to give an estimate of the quantity $\|\rho(t)/r\|_{L^2}$ which is very crucial in the proof of Theorem \ref{thm1}. The growth that we  will establish is  quadratic and this    improves the exponential growth  that one can easily obtain by writing the equation \mbox{of $\rho/r$.} More precisely we have
\begin{cor}\label{cor1}
Let $v$ be a smooth axisymmetric vector field  with zero divergence, $\rho_0\in L^2\cap L^\infty$ and $\rho$ be a solution of the transport equation
\begin{displaymath}
\left\{ \begin{array}{ll}
\partial_{t}\rho+v\cdot\nabla \rho=0&\\
{\rho}_{| t=0}=\rho_{0}.
\end{array} \right.
\end{displaymath}
Assume in addition that
$$d(\textnormal{supp }\rho_0,(Oz)):=r_0>0 \quad\hbox{and}\quad
\textnormal{diam}(\Pi_z(\textnormal{supp }\rho_0)):=d_0<\infty.$$ Then we have
$$
\int_{\RR^3}\frac{{\rho^2(t,x)}}{r^2}dx
\le
\frac{1}{r_0^2}\|\rho^0\|_{L^2}^2
+2\pi\|\rho_0\|_{L^\infty}^2\int_0^t\|({v^r}/{r})(\tau)\|_{L^\infty}d\tau
\Big(d_0+2\int_0^t\|v(\tau\|_{L^\infty}d\tau\Big),
$$
with $r=(x_1^2+x_2^2)^{\frac12}.$
\end{cor}
\begin{proof}
We have from the definition
\begin{eqnarray*}
\|(\rho/r)(t)\|_{L^2}^2
&=&
\int_{r\geq r_0}\frac{{\rho^2(t,x)}}{r^2}dx+\int_{r\le r_0}\frac{{\rho^2(t,x)}}{r^2}dx\\
&\le&
\frac{1}{r_0^2}\|\rho(t)\|_{L^2}^2+\|\rho(t)\|_{L^\infty}^2\int_{\{r\le r_0\}\cap
\textnormal{supp }\rho(t)}\frac{1}{r^2}dx\\
&\le&\frac{1}{r_0^2}\|\rho_0\|_{L^2}^2+\|\rho_0\|_{L^\infty}^2\int_{\{r\le r_0\}\cap
\textnormal{supp }\rho(t)}\frac{1}{r^2}dx.
\end{eqnarray*}
We have used the conservation of the $L^\infty$-norm of $\rho$.
Now using Proposition \ref{prop2} we get
\begin{eqnarray*}
\int_{\{r\le r_0\}\cap \textnormal{supp }\rho(t)}\frac{1}{r^2}dx
&\le&2\pi
\Big(\int_{r_0e^{-\int_0^t\|\frac{v^r}{r}(\tau)\|_{L^\infty}d\tau}
\le
r\le r_0}\frac{1}{r}dr\Big)\Big(\int_{\Pi_z(\textnormal{supp } \rho(t))}dz\Big)\\
&\le&2\pi\int_0^t\|\frac{v^r}{r}(\tau)\|_{L^\infty}d\tau\Big(d_0+2\int_0^t\|v(\tau\|_{L^\infty}d\tau\Big).
\end{eqnarray*}
This concludes the proof.
\end{proof}
%%%%%%%%%%%%%%
%% A priori estimates
\section{A priori estimates}
This section is devoted to the a priori estimates needed for the proof of Theorem \ref{thm1}. We distinguish especially two kinds: the first one deals with some easy estimates that one can obtained by energy estimates. However the second one is concerned with some strong estimates which are the heart of the proof of our main result.

\subsection{Weak a priori estimates}
We will prove the following energy estimates.
\begin{prop}\label{energy}
Let $v_0\in L^2$ be a  vector field  with zero divergence and $\rho_0\in L^2\cap L^\infty.$
Then every smooth solution of \eqref{bs} satisfies
$$
\|\rho(t)\|_{L^2\cap L^\infty}\le\|\rho_0\|_{L^2\cap L^\infty}.
$$
$$
\|v(t)\|_{L^2}^2+\int_0^t\|\nabla v(\tau)\|_{L^2}^2d\tau\le C_0(1+t^2).
$$

\end{prop}
\begin{proof}
The first estimate is obvious since the flow preserves Lebesgue measure.
For the second one we take the $L^2$-inner product of  the velocity equation with $v$.
Then we get after some integration by parts
\begin{equation}\label{eqs1}
\frac12\frac{d}{dt}\|v(t)\|_{L^2}^2+\|\nabla v(t)\|_{L^2}^2\le\|v(t)\|_{L^2}\|\rho(t)\|_{L^2}.
\end{equation}
After  simplification that one can rigorously  justify, the last inequality  leads   to
$$
\frac{d}{dt}\|v(t)\|_{L^2}\le\|\rho(t)\|_{L^2}.
$$
Integrating in time this inequality yields
$$
\|v(t)\|_{L^2}\le\|v_0\|_{L^2}+\int_0^t\|\rho(\tau)\|_{L^2}d\tau.
$$
Since  $\|\rho(t)\|_{L^2}=\|\rho_0\|_{L^2},$ then
$$
\|v(t)\|_{L^2}\le\|v_0\|_{L^2}+t\|\rho_0\|_{L^2}.
$$
Putting this estimate into \eqref{eqs1} gives
$$
\frac12\|v(t)\|_{L^2}^2+\int_0^t\|\nabla v(\tau)\|_{L^2}^2d\tau
\le
\frac12\|v_0\|_{L^2}^2+\big(\|v_0\|_{L^2}+t\|\rho_0\|_{L^2}\big)\|\rho_0\|_{L^2} t.
$$
This gives the desired estimate and the demonstration of the proposition is now accomplished.
\end{proof}
The next proposition describes some  estimates linking  the velocity to the
vorticity by the use of the so-called Bio-Savart law.
\begin{prop}\label{biot-s}
Let $v$ be a smooth  axisymmetric vector field with zero divergence  and denote
$\omega=\omega_\theta e_\theta$ its curl. Then
\begin{enumerate}
\item
$$
\|v\|_{L^\infty}\le C\|\omega_\theta\|_{L^{2}}^{\frac12}\|\omega_\theta\|_{\dot H^1}^{\frac12}.
$$
\item
$$
\|v^r/r\|_{L^\infty}\le C\|\omega_\theta/r\|_{L^{2}}^{\frac12}\|\omega_\theta/r\|_{\dot H^1}^{\frac12}.
$$
\end{enumerate}
\end{prop}
\begin{proof}
To start with, recall the classical  Biot-Savart law
$$
v(x)=\frac{1}{4\pi}\int_{\RR^3}\frac{(y-x)\wedge\omega(y)}{|y-x|^3}dy.
$$
It follows that
\begin{eqnarray*}
|v(x)|&\le&\frac{1}{4\pi}\big(\frac{1}{|\cdot|^2}\star|\omega|\big)(x)\\
&\lesssim&
 \frac{1}{|\cdot|^2}\star|\omega_\theta|:=J(x).
\end{eqnarray*}
Let $\lambda>0$ be a real number that will be fixed later.  Decompose the
convolution integral into two parts  as follows
\begin{eqnarray*}
J(x)&=&\int_{|x-y|\le\lambda}\frac{|\omega_\theta(y)|}{|x-y|^2}dy
+\int_{|x-y|\geq\lambda}\frac{|\omega_\theta(y)|}{|x-y|^2}dy\\
&:=&J_1(x)+J_2(x).
\end{eqnarray*}
Now using H\"{o}lder inequalities we get
\begin{eqnarray*}
\|J_1\|_{L^\infty}
&\le&
\|\omega_\theta\|_{L^6}\Big(\int_{|x|\le\lambda}\frac{1}{|x|^{\frac{12}{5}}}dx\Big)^{\frac56}\\
&\lesssim&\|\omega_\theta\|_{L^6}\lambda^{\frac12}.
\end{eqnarray*}
For the second integral we use again H\"{o}lder inequalities
\begin{eqnarray*}
\|J_2\|_{L^\infty}
&\le&
\|\omega_\theta\|_{L^2}\Big(\int_{|x|\geq\lambda}\frac{1}{|x|^{4}}dx\Big)^{\frac12}\\
&\lesssim&\|\omega_\theta\|_{L^2}\lambda^{-\frac12}.
\end{eqnarray*}
This yields
$$
\|v\|_{L^\infty}
\lesssim\
|\omega_\theta\|_{L^6}\lambda^{\frac12}+\|\omega_\theta\|_{L^2}\lambda^{-\frac12}.
$$
By choosing  $\lambda=\frac{\|\omega_\theta\|_{L^2}}{\|\omega_\theta\|_{L^6}}$ we get
$$
\|v\|_{L^\infty}\lesssim\|\omega_\theta\|_{L^2}^{\frac12}\|\omega_\theta\|_{L^6}^{\frac12}.
$$
It suffices now to use Sobolev embedding $\dot H^1(\RR^3)\hookrightarrow L^6$ to get the desired result.

For the second estimate we  recall the   following inequality, see for instance \cite{A-H-K,Taira},
$$
|v^r/r|(x)\lesssim \big( \frac{1}{|\cdot|^2}\star|\omega_\theta/r|\big)(x).
$$
Now we conclude similarly to the first estimate.
\end{proof}

\subsection{Strong a priori estimates}
The task is now  to find  some global estimates about more strong regularities of the solutions of \eqref{bs}.
The estimates developed below  will be the basic ingredient  of the proof of Theorem \ref{thm1}.
\begin{prop}\label{prop-cruc}
Let $v_0\in H^1$  be an axisymmetric vector field with zero
divergence such that $\frac{\textnormal{curl }{v_0}}{r}\in L^2.$ Let
$\rho_0\in L^2\cap L^\infty$ depending only on $(r,z)$ such that
$\hbox{supp }\rho_0$ does not intersect the axis $(Oz)$ and
$\Pi_z(\hbox{supp }\rho_0)$ is a compact set. Then every smooth
solution $(v,\rho)$ of the system \eqref{bs} satisfies for every
$t\geq 0$,
$$
\|v(t)\|_{H^1}^2+\int_0^t\|v(\tau)\|_{H^2}^2d\tau\le C_0e^{\exp{C_0t^9}},
$$
$$
\big\|\frac{\omega}{r}(t)\big\|_{L^2}^2+\big\|\frac{\omega}{r}\big\|_{L^2_t\dot H^1}^2
\le C_0 e^{\exp{C_0 t^9}}.
$$
Moreover, for every $p\in]3,\infty]$ we  have

\begin{equation}\label{eq15}
\| v\|_{L^1_t B_{p,1}^{\frac3p+1}}+\| \nabla v\|_{L^1_t L^\infty}\le  C_0 e^{\exp{C_0 t^9}}.
\end{equation}
The constant $C_0$ depends on the initial data.
\end{prop}
\begin{proof}
Taking the $L^2$-inner product of   the equation \eqref{tourbillon1} with $\omega_\theta$ we get

\begin{align*}
\frac12\frac{d}{dt}\|\omega_\theta\|_{L^2}^2+\|\nabla\omega_\theta\|_{L^2}^2
+\|{\omega_\theta\over r}\|_{L^2}^2
&=\int v^r{\omega_\theta\over r}\omega_\theta dx
-\int\partial_r\rho\,\omega_\theta dx.
\end{align*}
For the first integral  term of the r.h.s we  use  H\"{o}lder inequalities
$$
\int v^r{\omega_\theta\over r}\omega_\theta dx\le\|v\|_{L^6}\|{\omega_\theta\over r}\|_{L^2}
\|\omega_\theta\|_{L^3}.
$$
For the second term we integrate by parts taking into account that the density $\rho$
vanishes on the axis $(Oz)$ accordingly  to Corollary \eqref{cor1}
\begin{eqnarray*}
-\int\partial_r\rho\,\omega_\theta dx&=&-2\pi\int\partial_r\rho\,\omega_\theta rdr\,dz\\
&=&2\pi\int\rho\partial_r\omega_\theta\, rdr\,dz+2\pi\int\rho{\omega_\theta\over r} rdr\,dz\\
&=&\int\rho\big(\partial_r\omega_\theta+{\omega_\theta\over r}\big) dx.
\end{eqnarray*}
It follows from H\"{o}lder inequalities
$$
-\int\partial_r\rho\,\omega_\theta dx\le\|\rho\|_{L^2}\big(\|{\omega_\theta\over r}\|_{L^2}
+\|\partial_r\omega_\theta\|_{L^2}\big).
$$
Putting together these estimates yields
\begin{eqnarray*}
\frac12\frac{d}{dt}\|\omega_\theta\|_{L^2}^2+\|\nabla\omega_\theta\|_{L^2}^2
+\big\|{\omega_\theta\over r}\big\|_{L^2}^2
&\le&
\|v\|_{L^6}\big\|{\omega_\theta\over r}\big\|_{L^2}\|\omega_\theta\|_{L^3}
\\&
+&\|\rho\|_{L^2}\Big(\big\|{\omega_\theta\over r}\big\|_{L^2}+\|\partial_r\omega_\theta\|_{L^2}\Big).
\end{eqnarray*}
Using the inequality $|ab|\le \frac12 a^2+\frac{1}{2}b^2$ for the last term we obtain
\begin{eqnarray*}
\|\rho\|_{L^2}\big(\big\|{\omega_\theta\over r}\big\|_{L^2}
+\|\partial_r\omega_\theta\|_{L^2}\big)
&\le&
\|\rho\|_{L^2}^2+\frac12\big\|{\omega_\theta\over r}\big\|_{L^2}
^2+\frac12\|\partial_r\omega_\theta\|_{L^2}^2\\
&\le&\|\rho\|_{L^2}^2+\frac12\big\|{\omega_\theta\over r}\big\|_{L^2}
^2+\frac12\|\nabla\omega_\theta\|_{L^2}^2.
\end{eqnarray*}
We  have used in the last line the fact that in the cylindrical coordinates
$$
\|\nabla\omega_\theta\|_{L^2}^2=\|\partial_r\omega_\theta\|_{L^2}^2+\|\partial_z\omega_\theta\|_{L^2}^2.
$$
Therefore we get
\begin{equation}\label{eq05}
\frac{d}{dt}\|\omega_\theta\|_{L^2}^2+\|\nabla\omega_\theta\|_{L^2}^2
+\big\|{\omega_\theta\over r}\big\|_{L^2}^2
\le
2\|v\|_{L^6}\big\|{\omega_\theta\over r}\big\|_{L^2}\|\omega_\theta\|_{L^3}
+2\|\rho\|_{L^2}^2.
\end{equation}
Combining  the interpolation estimate
$\|\omega_\theta\|_{L^3}\lesssim\|\omega_\theta\|_{L^2}^{\frac12}
\|\nabla\omega_\theta\|_{L^2}^{\frac12}$
with Young inequality
$|ab|\le C_\eta
a^{\frac{1}{1-\eta}}+\frac{1}{4}b^{\frac1\eta},\,\forall
\eta\in]0,1[$ leads to
\begin{eqnarray}\label{eq06}
\nonumber\|v\|_{L^6}\big\|{\omega_\theta\over r}\big\|_{L^2}
\|\omega_\theta\|_{L^3}&\le &C\|v\|_{L^6}^{\frac43}\|{\omega_\theta\over r}\|_{L^2}^{\frac43}
\|\omega_\theta\|_{L^2}^{\frac23}+\frac14\|\nabla\omega_\theta\|_{L^2}^2
\\
&\le&
C\|\omega_\theta\|_{L^2}^2+\frac14\|\nabla\omega_\theta\|_{L^2}^2+C\|v\|_{L^6}^{2}
\big\|{\omega_\theta\over r}\big\|_{L^2}^{2}.
 \end{eqnarray}
 Inserting \eqref{eq06} into \eqref{eq05} and using the estimate
 $\|\rho(t)\|_{L^2}\le\|\rho_0\|_{L^2}$
\begin{equation*}
\frac{d}{dt}\|\omega_\theta\|_{L^2}^2+\frac12\|\nabla\omega_\theta\|_{L^2}^2
+\big\|{\omega_\theta\over r}\big\|_{L^2}^2
\le C\|\omega_\theta\|_{L^2}^2+C\|v\|_{L^6}^{2}\big\|{\omega_\theta\over r}\big\|_{L^2}^{2}
+2\|\rho_0\|_{L^2}^2.
\end{equation*}
Using Gronwall inequality and Proposition \ref{energy}  we get
\begin{eqnarray}\label{eq07}
\nonumber\|\omega_\theta(t)\|_{L^2}^2+\|\nabla\omega_\theta\|_{L^2_tL^2}^2
+\big\|{\omega_\theta\over r}\big\|_{L^2_tL^2}^2
\nonumber&\le& C e^{C t}\Big(\|\rho_0\|_{L^2}^2+\int_0^t\|v(\tau)\|_{L^6}^{2}
\big\|{\omega_\theta\over r}(\tau)\big\|_{L^2}^{2}d\tau\Big)\\
\nonumber&\le&C e^{C t}\Big(\|\rho_0\|_{L^2}^2+\|v\|_{L^2_t\dot H^1}^{2}
\big\|{\omega_\theta\over r}\big\|_{L^\infty_tL^2}^{2}\Big)\\
&\le&C_0e^{Ct}\Big(1+\big\|{\omega_\theta\over r}\big\|_{L^\infty_tL^2}^{2}\Big).
\end{eqnarray}
Let us now move to the equation \eqref{equation_i} and try to estimate in a proper way the \mbox{quantity $\Gamma:=\frac{\omega_\theta}{r}$.}
Then taking the $L^2$-inner product of  \eqref{equation_i} with $\Gamma$  and integrating by parts using the incompressibility of the
flow and $\rho(t,0,z)=0$ gives
\begin{align*}
\frac12\frac{d}{dt}\|\Gamma\|_{L^2}^2
+\|\partial_r\Gamma\|_{L^2}^2
+\|\partial_z\Gamma\|_{L^2}^2
-4\pi\int\partial_r(\Gamma)
{\Gamma} drdz
&=
-2\pi\int\partial_r\rho\,\Gamma drdz
\\&
=2\pi\int{\rho\over r}\,\partial_r\Gamma\, r drdz
\\&
\le
\big\|{\rho\over r}\big\|_{L^2}\| \partial_r\Gamma\|_{L^2}.
\end{align*}
Since
$$
4\pi\int\partial_r(\Gamma)\Gamma drdz
=2\pi\int_{\RR}\int_0^{+\infty}\partial_r(\Gamma)^2 drdz
\le
0,
$$
then and once again by Young inequality yields
\begin{equation*}
\frac{d}{dt}\big\|{\omega_\theta\over r}\big\|_{L^2}^2
+\big\|\nabla\big({\omega_\theta\over r}\big)\big\|_{L^2}^2
\le
\big\|{\rho\over r}\big\|_{L^2}^2.
\end{equation*}
Integrating in time this differential inequality we get
\begin{equation}\label{eq08}
\big\|{\omega_\theta\over r}(t)\big\|_{L^2}^2
+\int_0^t\big\|\nabla\big({\omega_\theta\over r}\big)(\tau)\big\|_{L^2}^2d\tau
\le\big\|{\omega_\theta\over r}(0)\big\|_{L^2}^2+
\int_0^t\big\|{\rho\over r}(\tau)\big\|_{L^2}^2d\tau.
\end{equation}
To estimate the last term of the above inequality we will use Corollary \ref{cor1}
$$
\big\|\frac{\rho}{r}(t)\big\|_{L^2}^2
\le
C_0+C_0\int_0^t\big\|\frac{v^r}{r}(\tau)\big\|_{L^\infty}d\tau
\Big(1+\int_0^t\|v(\tau)\|_{L^\infty}d\tau\Big).
$$
Thus
\begin{align*}
\big\|{\omega_\theta\over r}(t)\big\|_{L^2}^2
&+\big\|\nabla\big({\omega_\theta\over r}\big)\big\|_{L^2_tL^2}^2
\le
C_0 (1+t)
\\&
+C_0\int_0^t\Big\{\int_0^{t'}\|(\frac{v^r}{r})(\tau)\|_{L^\infty}d\tau
\Big(1+\int_0^{t'}\|v(\tau)\|_{L^\infty}d\tau\Big)\Big\}dt'
\\&
\le
C_0(1+t)+C_0t\int_0^t\|(\frac{v^r}{r})(\tau)\|_{L^\infty}d\tau
\\&+C_0 t\int_0^t\|(\frac{v^r}{r})(\tau)\|_{L^\infty}d\tau\,\int_0^{t}\|v(\tau)\|_{L^\infty}d\tau.
\end{align*}
From Young and H\"{o}lder inequalities we get successively
\begin{eqnarray*}
C_0t\int_0^t\|(\frac{v^r}{r})(\tau)\|_{L^\infty}d\tau
&\le &\int_0^t\big(C_0 t^2+\|(\frac{v^r}{r})(\tau)\|_{L^\infty}^2\big)d\tau\\
&\le&C_0 t^3+\int_0^t\|(\frac{v^r}{r})(\tau)\|_{L^\infty}^2d\tau
\end{eqnarray*}
and
\begin{eqnarray*}
C_0 t\|{v^r}/{r}\|_{L^1_tL^\infty}\|v\|_{L^1_tL^\infty}
&\le&
C_0 t^2\|v^r/r\|_{L^2_tL^\infty}\|v\|_{L^2_tL^\infty}\\
&\le&
C_0 t^4\int_0^t \|(\frac{v^r}{r})(\tau)\|_{L^\infty}^2d\tau+\int_0^t \|v(\tau)\|_{L^\infty}^2d\tau.
\end{eqnarray*}
Putting together these estimates yields
\begin{align}\label{eq10}
\nonumber\big\|{\omega_\theta\over r}(t)\big\|_{L^2}^2
&+\big\|\nabla\big({\omega_\theta\over r}\big)\big\|_{L^2_tL^2}^2
\le
C_0 (1+t^4)\Big(1+\int_0^t\big\|({v^r}/{r})(\tau)\big\|_{L^\infty}^2d\tau\Big)
\\&
+\int_0^{t}\|v(\tau)\|_{L^\infty}^2d\tau.
\end{align}
Recall from Proposition \ref{biot-s} that
\begin{eqnarray*}
\big\|{v^r}/{r}\big\|_{L^\infty}
&\lesssim&
\big\|\frac{\omega_\theta}{r}\big\|_{L^2}^{\frac12}
\big\|\nabla(\frac{\omega_\theta}{r})\big\|_{L^2}^{\frac12}.
\end{eqnarray*}
Accordingly,  we get by Young inequalities
\begin{equation*}
C_0(1+t^4)\int_0^{t}\|({v^r}/{r})(\tau)\|_{L^\infty}^2d\tau
\leq
C_0(1+t^8)\int_0^{t}\big\|\frac{\omega_\theta}{r}(\tau)\big\|_{L^2}^2d\tau
+\frac12\int_0^{t}\big\|\nabla(\frac{\omega_\theta}{r})(\tau)\big\|_{L^2}^2d\tau.
\end{equation*}
Inserting this estimate into \eqref{eq10}
\begin{align}\label{eq100}
\nonumber\big\|{\omega_\theta\over r}(t)\big\|_{L^2}^2
+\big\|\nabla\big({\omega_\theta\over r}\big)\big\|_{L^2_tL^2}^2
&\le
C_0(1+t^8)\int_0^{t}\big\|\frac{\omega_\theta}{r}(\tau)\big\|_{L^2}^2d\tau
\\&
+\int_0^{t}\|v(\tau)\|_{L^\infty}^2d\tau.
\end{align}

Now using Gronwall inequality we find
\begin{equation}\label{eq11}
 \big\|{\omega_\theta\over r}(t)\big\|_{L^2}^2
+\big\|\nabla\big({\omega_\theta\over r}\big)\big\|_{L^2_tL^2}^2
\le
C_0e^{C_0t^9}\Big(1+\int_0^{t}\|v(\tau)\|_{L^\infty}^2d\tau\Big).
\end{equation}
Putting\eqref{eq11} into \eqref{eq07} yields
\begin{equation}\label{eq12}
\|\omega_\theta(t)\|_{L^2}^2+\|\nabla\omega_\theta\|_{L^2_tL^2}^2
+\big\|{\omega_\theta\over r}\big\|_{L^2_tL^2}^2
\le C_0 e^{C_0 t^9}\Big(1+\int_0^{t}\|v(\tau)\|_{L^\infty}^2d\tau\Big).
\end{equation}
To estimate the last term in the r.h.s of (\ref{eq12}) we use Proposition \ref{biot-s}
\begin{eqnarray*}
\|v\|_{L^\infty}
&\lesssim&
\|\omega_\theta\|_{L^2}^{\frac12}\|\nabla\omega_\theta\|_{L^2}^{\frac12}.
\end{eqnarray*}
It follows
\begin{equation}\label{born}
\int_0^t\|v(\tau\|_{L^\infty}^2d\tau\le\int_0^t\|\omega_\theta(\tau\|_{L^2}
\|\nabla\omega_\theta(\tau)\|_{L^2}d\tau.
\end{equation}
Therefore we get from \eqref{eq12} and Young inequality
\begin{eqnarray*}
\|\omega_\theta(t)\|_{L^2}^2+\|\nabla\omega_\theta\|_{L^2_tL^2}^2
+\big\|{\omega_\theta\over r}\big\|_{L^2_tL^2}^2 &\le& C_0 e^{C_0
t^9}\Big(1+\int_0^t\|\omega_\theta(\tau\|_{L^2}\|\nabla\omega_\theta(\tau)\|_{L^2}d\tau\Big)
\\
&\le& C_0 e^{C_0 t^9}\Big(1+\int_0^t\|\omega_\theta(\tau\|_{L^2}^2d\tau\Big)
+\frac12\|\nabla\omega_\theta\|_{L^2_tL^2}^2.
\end{eqnarray*}
It suffices now to use Gronwall inequality
$$
\|\omega_\theta(t)\|_{L^2}^2+\|\nabla\omega_\theta\|_{L^2_tL^2}^2
+\big\|{\omega_\theta\over r}\big\|_{L^2_tL^2}^2
\le C_0 e^{\exp{C_0 t^9}}.
$$
Since in cylindrical coordinates we have $\|\nabla\omega\|_{L^2}^2=
\|\nabla\omega_\theta\|_{L^2}^2+\big\|\frac{\omega_\theta}{r}\big\|_{L^2}^2$
then
\begin{equation}\label{born2}
\|\omega(t)\|_{L^2}^2+\|\nabla\omega\|_{L^2_tL^2}^2\le C_0 e^{\exp{C_0 t^9}}.
\end{equation}
Inserting this estimate into  \eqref{born} gives
\begin{equation}\label{born1}
\|v\|_{L^2_tL^\infty}\le C_0e^{\exp{C_0t^9}}.
\end{equation}
Combining this estimate with \eqref{eq11} and \eqref{born2} yields
$$
\big\|\frac{\omega_\theta}{r}(t)\big\|_{L^2}^2
+\big\|\nabla\big(\frac{\omega_\theta}{r}\big)(t)\big\|_{L^2_tL^2}^2
\le C_0 e^{\exp{C_0 t^9}}.
$$
This concludes the first part of Proposition \ref{prop-cruc}.
Let us now show how to prove  the \mbox{estimate \eqref{eq15}.} Let $q\in\NN$ and set $v_q:=\Delta_qv$.  Then localizing in frequency the first
equation of  \eqref{bs} and using Duhamel formula we get
$$
v_q(t)=e^{t\Delta}v_q(0)+\int_{0}^te^{(t-\tau)\Delta}
\Delta_q\big(\mathcal{P}(v\cdot \nabla v)(\tau)d\tau
+\int_{0}^te^{(t-\tau)\Delta}\Delta_q\big(\mathcal{P}(\rho e_z))(\tau)d\tau,
$$
where $\mathcal{P}$ denotes Leray's projector over solenoidal vector fields. Now we will use a local version of the  smoothing effects of the heat semigroup, for the proof see for example \cite{che1},
$$
\|e^{t\Delta}\Delta_q\mathcal{P}f\|_{L^p}\le Ce^{-ct2^{2q}}
\|\Delta_qf\|_{L^p},\quad\forall p\in[1,\infty].
$$
Combining this estimate with Bernstein inequality gives
\begin{eqnarray*}
\|v_q(t)\|_{L^p}\lesssim e^{-ct2^{2q}}\|v_q(0)\|_{L^p}
&+&
2^q\int_0^t e^{-c(t-\tau)2^{2q}}\|\Delta_q(v\otimes v)(\tau)\|_{L^p}d\tau
\\&+&
\int_0^t e^{-c(t-\tau)2^{2q}}\|\Delta_q\rho(\tau)\|_{L^p}d\tau.
\end{eqnarray*}
Integrating in time and using convolution inequalities  lead to
$$
\|v_q\|_{L^1_tL^p}\lesssim 2^{-2q}\|v_q(0)\|_{L^p}
+2^{-q}\int_0^t\|\Delta_q(v\otimes v)(\tau)\|_{L^p}d\tau
+2^{-2q}\int_0^t\|\Delta_q\rho(\tau)\|_{L^p}d\tau.
$$
It follows that
\begin{eqnarray*}
\|v\|_{L^1_tB_{p,1}^{\frac3p+1}}
\le
\|\Delta_{-1}v\|_{L^1_tL^p}
+\|v_0\|_{B_{p,1}^{\frac3p-1}}&+&\int_0^t\|(v\otimes v)(\tau)\|_{B_{p,1}^{\frac3p}}d\tau
\\&+&
\int_0^t\|\rho(\tau)\|_{B_{p,1}^{\frac3p-1}}d\tau.
\end{eqnarray*}
For the first term of the r.h.s we combine Bernstein inequality with
Proposition \ref{energy} ($p\geq 2$)
\begin{eqnarray*}
\|\Delta_{-1}v\|_{L^1_tL^p}&\lesssim &t\|v\|_{L^\infty_t L^2}\\
&\le&C_0(1+t^2).
\end{eqnarray*}
For the last term of the r.h.s we use the embedding
$L^p\hookrightarrow B_{p,1}^{\frac3p-1},\,$ for $p>3$, combined with Proposition \ref{energy}
\begin{eqnarray*}
\int_0^t\|\rho(\tau)\|_{B_{p,1}^{\frac3p-1}}d\tau&\lesssim& t\|\rho\|_{L^\infty_tL^p}\\
&\le& C_0 t.
\end{eqnarray*}
On the other hand we have from Besov embedding
\begin{eqnarray*}
\|v_0\|_{B_{p,1}^{\frac3p-1}}&\lesssim& \|v_0\|_{B_{2,1}^{\frac12}}\\
&\le& \|v_0\|_{H^1}.
\end{eqnarray*}
Putting together these inequalities we find
$$
\|v\|_{L^1_tB_{p,1}^{\frac3p+1}}\le C_0(1+t^2)+\int_0^t\|(v\otimes v)(\tau)\|_{B_{p,1}^{\frac3p}}d\tau.
$$
Now we use Besov embeddings, law products and interpolation results
\begin{eqnarray*}
\|v\otimes v\|_{B_{p,1}^{\frac3p}}&\lesssim&\|v\otimes v\|_{B_{2,1}^{\frac32}}\\
&\lesssim&\| v\|_{L^\infty}\|v\|_{B_{2,1}^{\frac32}}\\
&\lesssim&\| v\|_{L^\infty}\|v\|_{L^2}+\| v\|_{L^\infty}\|\omega\|_{B_{2,1}^{\frac12}}
\\
&\lesssim&
\| v\|_{L^\infty}\|v\|_{L^2}+\| v\|_{L^\infty}\|\omega\|_{L^2}^{\frac12}\|\omega\|_{H^1}^{\frac12}.
\end{eqnarray*}
This yields according to H\"{o}lder inequalities
\begin{eqnarray*}
\|v\otimes v\|_{L^1_tB_{p,1}^{\frac3p}}
&\lesssim&
t^{\frac12}\|v\|_{L^2_tL^\infty}
\|v\|_{L^\infty_t L^2}+\|v\|_{L^{\frac43}_tL^\infty}\|\omega\|_{L^\infty_tL^2}^{\frac12}
\|\omega\|_{L^2_tH^1}^{\frac12}\\
&\lesssim&
t^{\frac12}\|v\|_{L^2_tL^\infty}\|v\|_{L^\infty_t L^2}
+t^{\frac14}\|v\|_{L^{2}_tL^\infty}\|\omega\|_{L^\infty_tL^2}^{\frac12}
\|\omega\|_{L^2_tH^1}^{\frac12}.
\end{eqnarray*}
It suffices now to use \eqref{born2} and \eqref{born1}
$$
\|v\otimes v\|_{L^1_tB_{p,1}^{\frac3p}}\le C_0e^{\exp{C_0 t^9}}.
$$
Hence we obtain
$$
\|v\|_{L^1_tB_{p,1}^{\frac3p+1}}\le C_0e^{\exp{C_0t^9}}.
$$
The last estimate of Proposition \eqref{eq15} is a direct consequence of the
\mbox{embedding $B_{p,1}^{\frac3p+1}\hookrightarrow W^{1,\infty}$} and the proof is then accomplished.
\end{proof}
\section{Proof of the main result}
The existence part can be done in a  classical way for example  by smoothing out the initial data.
However we have to use an  appropriate approximation that does not alter  never the initial
geometric structure nor the uniform estimates in the space of initial data.
We will work with the approximation of the identity and show that it has the
requested properties: let $\phi$  be a smooth positive radial function with
support contained in $B(0,1)$ and such that $\phi\equiv1$ in a neighborhood of zero.
We assume \mbox{that $\displaystyle{\int_{\RR^3}\phi(x)dx=1.}$} For every $n\in\NN^*$
we set $\phi_n(x)={n^3}\phi(nx)$ and we define   the family
$$
v_{0,n}=\phi_n \star v_0\quad\textnormal{and}\quad \rho_{0,n}=\phi_n\star\rho_0.
$$
We will start with  the following stability results.
\begin{lem}\label{lem1}
\begin{enumerate}
\item
Let $v_0\in H^1$ be an axisymmetric vector field  with zero divergence and such that
$({\textnormal{curl }v_0})/{r}\in L^2$. Then  for every $n\in \NN^*$ the vector field $v_{0,n}$ is
axisymmetric with zero divergence. Moreover, there exists a constant $C$ depending on $\phi$ such that
$$
\|v_{0,n}\|_{H^1}\le \|v_0\|_{H^1},\quad
\big\| (\textnormal{curl }v_{0,n})/r \big\|_{L^2}\le C\big\| {(\textnormal{curl }v_0)}/{r} \big\|_{L^2}
$$
\item Let $v$ be a smooth axisymmetric vector field with zero divergence  and $\rho$ be
a solution of the transport equation
\begin{displaymath}
\left\{ \begin{array}{ll}
\partial_{t}\rho^n+v\cdot\nabla \rho^n =0&\\
{\rho^n}_{| t=0}=\rho_{0,n}.
\end{array} \right.
\end{displaymath}
Assume that $d(\textnormal{supp }\rho_0,(Oz)):=r_0>0$ and
${\textnormal{diam }}(\Pi_z(\textnormal{supp }\rho_0)):=d_0<\infty$.
Then there exists constants  $n_0$  and $C$ depending only on $r_0, d_0$ and $\phi$ such that
$$
\int_{\RR^3}\frac{{{(\rho^n)}^2(t,x)}}{r^2}dx
\le
\frac{1}{r_0^2}\|\rho_0\|_{L^2}^2+2\|\rho_0\|_{L^\infty}
\big(\ln 2+\|{v^r}/{r}\|_{L^1_tL^\infty}\big)\big(d_0+\|v\|_{L^1_tL^\infty}\big).
$$

\end{enumerate}
\end{lem}
\begin{Rema}
The estimate  of the part (2) of the above lemma is  little bit different from part (2)
of Corollary \ref{cor1} because we have an   additional  linear term $\|v\|_{L^1_t L^\infty}$.
Nevertheless, the presence of this term does not  deeply affect the calculus seen before in the a priori estimates;  we obtain at  the end
the same estimates of Proposition \ref{prop-cruc}.
\end{Rema}
\begin{proof}
{\bf (1)} The  fact that the vector field $v_{0,n}$ is axisymmetric is due to the radial property  of the functions
$\phi_n$, for more details see \cite{A-H-K}. The estimate of $v_{0,n}$ in  $H^1$ is
easy to obtain by using the classical properties of the convolution  operation combined with
$\|\phi_n\|_{L^1}=1.$ Concerning the demonstration of  the second estimate, it is more subtle and  we refer to \cite{benamer} where it
is proven for more general framework Lebesgue space, that is $L^p$, \mbox{for all $p\in[1,\infty]$.}
\

{\bf(2)} In order to establish  the desired   estimate we need to check  that the estimates of Proposition \ref{prop2}
are stable with respect to $n$. More precisely, we will prove that  for sufficiently large $n\geq n_0$
\begin{equation}\label{eq765}
d (\textnormal{supp }\rho^n(t),(Oz))\geq\, \frac{r_0}{2}\,
e^{-\int_0^t\|\frac{v^r}{r}(\tau)\|_{L^\infty}d\tau}
\end{equation}
and
\begin{equation}\label{eq764}
d_n(t)\le 2d_0+2\int_0^t\|v(\tau\|_{L^\infty}d\tau\quad
{\textnormal{ with}}\quad d_n(t):={\textnormal{diam }}(\Pi_z(\textnormal{supp }\rho^n (t))).
\end{equation}
Assume that we have the inequalities  \eqref{eq765} and  \eqref{eq764}
then reproducing the proof of Corollary \ref{cor1} we get
\begin{eqnarray*}
\|(\rho_n/r)(t)\|_{L^2}^2
&\le&
\frac{1}{r_0^2}\|\rho_{0,n}\|_{L^2}^2+\|\rho_{0,n}\|_{L^\infty}^2\int_{\{r\le r_0\}
\cap\textnormal{supp }\rho_n(t)}\frac{1}{r^2}dx
\\&\le&
\frac{1}{r_0^2}\|\rho_{0}\|_{L^2}^2+\|\rho_{0}\|_{L^\infty}^2\int_{\{r\le r_0\}
\cap\textnormal{supp }\rho_n(t)}\frac{1}{r^2}dx
\\&\le&
\frac{1}{r_0^2}\|\rho_{0}\|_{L^2}^2+\|\rho_{0}\|_{L^\infty}^2
\big(\ln(2)+\|v^r/r\|_{L^1_tL^\infty}\big)\big(2d_0+2\int_0^t\|v(\tau\|_{L^\infty}d\tau  \big).
\end{eqnarray*}
This is what we want. Now let us come back  to the proof of the inequality \eqref{eq765}.
Following the same proof of Proposition \ref{prop2} it appears that one needs only
to establish  the stability  conditions of the support, that is, for large $n\geq n_0,$

 \begin{equation}\label{eq777}
d(\textnormal{supp }\rho_0^n,(Oz))\geq\, \frac{r_0}{2}.
\end{equation}
From the definition we have
$$
\textnormal{supp }\rho_0^n\subset\Big\{x;\|x-y\|
\le
\frac1n \textnormal{ for some },  \|y'\|\geq r_0   \Big\},\quad
\textnormal{ with } y=(y',y_3), y'\in\RR^2.
$$
It is obvious that for $x\in \textnormal{supp }\rho_0^n$ we have
$\|x'-y'\|\le\frac1n$ and then by the triangular inequality
$$
\|x'\|\geq r_0-\frac1n\cdot
$$
Thus there exists $n_0\in \NN$ such that for every $n\geq n_0$ we have
$\|x'\|\geq \frac{r_0}{2}\cdot$ This proves that
$\textnormal{supp }\rho_0^n\subset\{x, \|x'\|\geq\frac{r_0}{2} \}$ and then we deduce \eqref{eq777}.

Let us now show how to get  the estimate \eqref{eq764}. According to part (2)
of \mbox{Proposition \ref{prop2}} we have
$$
d_n(t)\le d_{0,n}+2\int_0^t\|v(\tau)\|_{L^\infty}d\tau,
$$
with $d_{0,n}:=\textnormal{diam }(\Pi_z(\textnormal{supp }\rho_{0,n})).$
Now it is easy to see from the support properties of the convolution  that
$$
d_{0,n}\le d_0+\frac2n\cdot
$$
Taking $n$ large enough, $n\geq n_1$ for some $n_1$,  yields  $d_{0,n}\le 2d_0.$
Now to be sure that both inequalities  \eqref{eq765} and \eqref{eq764} are simultaneously
satisfied we take $n\geq \max\{n_0,n_1\}.$

\end{proof}
Let us now come back to the proof of the existence part of Theorem
\ref{thm1}. We have just seen in  Lemma \ref{lem1} that  the initial
structure of axisymmetry is preserved for every $n$ and the involved
norms are uniformly controlled with respect to this parameter $n$.
Thus we can construct locally in time a unique solution
$(v_n,\theta_n)$ that does not blow up in finite time since the
Lipschitz norm of the velocity is well controlled as it was stated
in \mbox{Proposition \ref{prop-cruc}.} By standard arguments we can
show that this family converges to $(v,\theta)$ which satisfies in
turn our IVP. We omit here the details which are very classical and we will next focus on the
uniqueness part. Set
$$
\mathcal{X}_T:=\big(L^\infty_TL^2\cap L^2_TH^1\big)\times L^\infty_T H^{-1}.
$$
Let $(v^i,\rho^i)\in \mathcal{X}_T, 1\leq i\leq 2$ be  two solutions of the
 system \eqref{bs} with the same initial \mbox{data $(v_0,\theta_0)$} and denote
 $\delta v=v^2-v^1,\delta\theta=\theta^2-\theta^1$. Then
\begin{equation}\label{diff}
\left\{ \begin{array}{ll}
\partial_t\delta v+v^2\cdot\nabla\delta v-\Delta\delta
v+\nabla\delta p
=-\delta v\cdot\nabla v^1+\delta\rho e_z
\\
\partial_t \delta\rho +v^2\cdot\nabla\delta\rho=
 -\delta v\cdot\nabla \rho^1
\\
{\mathop{\rm div}}\,  v^i=0.
\end{array}
\right.
\end{equation}
Taking the $L^2$-inner product of the first equation with $\delta v$  and integrating by parts
$$
\begin{aligned}
\frac{1}{2}\frac{d}{dt}\|\delta v\|_{L^2}^2
+\|\nabla\delta v\|_{L^2}^2
&=-\int\delta v\cdot\nabla v^1\delta v\,dx
+\int\delta\rho e_z\delta v\,dx
\\&
\le\|\delta v\|_{L^2}\|v^1\|_{L^\infty}\|\nabla\delta v\|_{L^2}
+\|\delta\rho\|_{H^{-1}}\|\delta v\|_{H^1}
\\&
\lesssim
\big(\|\delta v\|_{L^2}\|v^1\|_{L^\infty}
+\|\delta\rho\|_{H^{-1}}\big)
\|\nabla\delta v\|_{L^2}+\|\delta\rho\|_{H^{-1}}\|\delta v\|_{L^2}.
\end{aligned}
$$
Applying  Young inequality leads to
\begin{equation}\label{diff}
\frac{d}{dt}\|\delta v\|_{L^2}^2
+\|\nabla\delta v\|_{L^2}^2
\lesssim
\|\delta v\|_{L^2}^2\big(\|v^1\|_{L^\infty}^2+1\big)
+\|\delta\rho\|_{ H^{-1}}^2.
\end{equation}
To estimate $\|\delta\rho\|_{H^{-1}}$ we will use
Proposition 3.1 of  \cite{AP}: for every $p\in[2,\infty[$
$$
\|\delta\rho(t)\|_{H^{-1}}
\le
C\|\delta v\cdot\nabla\rho^1\|_{L^1_tH^{-1}}
\exp{(C\|\nabla v^2\|_{L^1_tB^{\frac{3}{p}}_{p,1}})}.
$$
As ${\mathop{\rm div}}\, \delta v=0,$ then
\begin{eqnarray*}
\|\delta v\cdot\nabla\rho^1\|_{L^1_tH^{-1}}&\le&\|\delta v\,\rho^1\|_{L^1_tL^2}\\
&\le& \|\rho_0\|_{L^\infty}\|\delta v\|_{L^1_tL^2}.
\end{eqnarray*}
Inserting this estimate into
(\ref{diff}), we obtain
\begin{eqnarray*}
\frac{d}{dt}\|\delta v(t)\|_{L^2}^2
+\|\nabla\delta v(t)\|_{L^2}^2
&\le&
C\|\delta v(t)\|_{L^2}^2\big(\|v^1(t)\|_{L^\infty}^2+1\big)\\
&+&C\exp{(C\|\nabla v^2\|_{L^1_tB^{\frac{3}{p}}_{p,1}})}
\|\rho_0\|_{L^\infty}^2\|\delta v\|_{L^1_tL^2}^2.
\end{eqnarray*}
Integrating this differential inequality we get
\begin{eqnarray*}
\|\delta v\|_{L^\infty_tL^2}^2
&\le&
C\int_0^t\|\delta v(\tau)\|_{L^2}^2\big(\|v^1(\tau)\|_{L^\infty}^2+1\big)d\tau
\\&+&
C\exp(C\|\nabla v^2\|_{L^1_tB^{\frac{3}{p}}_{p,1}})\|\rho_0\|_{L^\infty}^2
\int_0^t\|\delta v\|_{L^1_\tau L^2}^2d\tau
\\&
\le&C\int_0^t\|\delta v\|_{L^\infty_\tau L^2}^2\big(\|v^1(\tau)\|_{L^\infty}^2+1\big)d\tau
\\&+&
C\exp(C\|\nabla v^2\|_{L^1_tB^{\frac{3}{p}}_{p,1}})
\|\rho_0\|_{L^\infty}t^2\int_0^t\|\delta v\|_{L^\infty_\tau L^2}^2d\tau.
\end{eqnarray*}
It suffices now to use Gronwall inequality.

 \end{document}